\documentclass[11pt]{article}
\usepackage{amscd,amssymb,stmaryrd}
\usepackage[fleqn]{amsmath}
\usepackage{graphicx,graphics,epstopdf}
\usepackage{graphicx}
\graphicspath{{./Figure/Figure_partial/}}
\usepackage{float}
\usepackage{subcaption}
\usepackage[utf8]{inputenc}
\usepackage[export]{adjustbox}
\usepackage{wrapfig}
\usepackage{amstext}
\usepackage{pstricks,pst-node,pst-plot,pst-coil}
\usepackage{amsthm}
\usepackage{amsmath}
\usepackage{color}
\usepackage{mathtools}
\usepackage{hyperref}
\usepackage[english]{babel}
\usepackage{booktabs}
\usepackage{mathrsfs}
\usepackage{comment}
\usepackage{marginnote}
\usepackage{enumitem}

\newcommand\norm[1]{\left\|#1\right\|}

\theoremstyle{definition}
\newtheorem{theorem}{Theorem}[section]

\theoremstyle{remark}

\title{Partially Explicit Time Discretization for Time Fractional Diffusion Equation}
\author{Jiuhua Hu\footnote{Department of Mathematics, Texas A\&M University, College Station, TX 77843, USA. (\texttt{E-mail: jiuhuahu@tamu.edu})}, ~Anatoly Alikhanov\footnote{North-Caucasus Center for Mathematical Research, North-Caucasus Federal University,  Stavropol, 355017  Russia. (\texttt{E-mail: aaalikhanov@gmail.com})}, ~Yalchin Efendiev\footnote{Department of Mathematics, Texas A\&M University, College Station, TX 77843, USA \& North-Eastern Federal University, Yakutsk, Russia. (\texttt{E-mail: efendiev@math.tamu.edu})}, ~ Wing Tat Leung\footnote{Department of Mathematics, University of California, Irvine, USA. (\texttt{E-mail: wtleung@uci.edu})}}
\date{\today} 
\begin{document}
	\maketitle
	\begin{abstract}

Time fractional PDEs have been used in many applications for modeling and simulations. Many of these applications are multiscale and contain high contrast 
variations in the media properties.  It requires
very small time step size to perform detailed computations. On the other hand, in the presence of small spatial
grids, very small time step size is required for explicit methods. Explicit methods have many advantages as we discuss in the paper. 
In this paper, we propose a partial explicit method for time fractional PDEs. 
The approach solves the forward problem
on a coarse computational grid, which is much larger than 
spatial heterogeneities,  and requires only a few degrees of freedom
to be treated implicitly. Via the construction of appropriate spaces
and careful stability analysis, we can show that the time step can 
be chosen not to depend on the contrast or scale as the coarse mesh size.
Thus, one can use larger time step size in an explicit approach. 
We present stability theory for our proposed method and our numerical results confirm the stability findings and demonstrate the performance of the approach.

	\end{abstract}
	

	\section{Introduction}


Many problems have multiscale nature. These include flow in porous media,
composite materials and so on. 
Multiscale features typically occur spatially due to variable nature
of media properties. For example, porous media properties can vary
at different scales. Many multiscale methods have been developed
for steady state and dynamic problems. However, there has been 
limited research on multiscale methods with fractional time derivatives.
Fractional time derivatives occur in many applications
as will be discussed next.


Recently, there has been many research activities related to equations 
with fractional derivatives. This is due to the numerous effective 
applications of fractional calculus in various fields of science 
and technology \cite{kilbas2006theory,hilfer2000applications,oldham1974fractional,podlubny1998fractional}. For example, fractional 
derivatives are used for describing the physical process of 
statistical transfer and this leads to diffusion equations 
of fractional orders \cite{chukbar1994stochastic,nigmatullin1986realization}. Many methods have been developed for the 
numerical solution of time fractional equations. 
The most common 
differential approximation of the time fractional derivative is the so-called 
L1 method \cite{oldham1974fractional,zhang2014finite}. To approximate the time fractional derivative with a higher 
order of accuracy, difference analogs of the L2 type are used \cite{alikhanov2015new,alikhanov2021high,gao2014new}. 
The basic properties of these difference analogs are well studied
and are used for solving various fractional PDEs.
In this paper, we combine time fractional PDE approximation with 
multiscale splitting to solve challenging 
time fractional PDEs in heterogeneous media.

For time fractional PDEs, it is difficult to use explicit methods due to
time step constraints. The time step constraints involve the fractional power
and the fine grid which is much smaller compared to the coarse grid
and thus one requires smaller time step size for explicit 
methods as the fractional power gets smaller. In the presence of 
multiscale features, the time step size is even smaller as the contrast
increases.

 In this paper, we propose 
partially explicit approach, where one implicitly treats a few degrees of 
freedom defined on a coarse grid and the rest of degrees of freedom are
treated explicitly. This approach allows removing the time step constraint
due to the contrast and reduce the time step constraint due to the fine-grid
mesh size. In particular, the time step size is smaller than
an appropriate fractional power of coarse-grid mesh size. In general, 
by choosing coarse-grid mesh size larger, one can allieviate this problem
and make the time step larger. We would like to note that explicit methods
have many advantages as they provide less communications and easy to
compute. They can also be used in constructing efficient neural network
architectures.

In the paper, we present a novel framework for stability analysis
of implicit-explicit methods for time fractional PDEs. 
The stability analysis starts with a space decomposition for 
coarse-grid degrees of freedom and the other degrees of freedom.
This framework 
derives conditions necessary for partial explicit schemes to be stable.
In particular, the stability conditions show that one needs
the second space to be free of contrast. To achieve this, we 
need spaces like CEM-GMsFEM, which we describe next.

In previous findings, many multiscale algorithms have been developed, 
such as
RVE based homogenization approaches \cite{eh09,le2014msfem}, 
multiscale
finite element methods \cite{eh09,hw97,jennylt03}, 
generalized multiscale finite element methods (GMsFEM) \cite{chung2016adaptiveJCP,MixedGMsFEM,WaveGMsFEM,chung2018fast,GMsFEM13}, 
constraint energy minimizing GMsFEM (CEM-GMsFEM) 
\cite{chung2018constraint, chung2018constraintmixed}, nonlocal
multi-continua (NLMC) approaches \cite{NLMC},
metric-based upscaling \cite{oz06_1}, heterogeneous multiscale method 
\cite{ee03}, localized orthogonal decomposition (LOD) 
\cite{henning2012localized}, equation-free approaches \cite{rk07,skr06}, 
multiscale stochastic approaches \cite{hou2018adaptive,hou2017exploring, hou2019model},
and hierarchical multiscale method \cite{brown2013efficient},
are developed to address spatial heterogeneities.
For high-contrast problems, approaches that require multiple
multiscale basis functions are needed, which include
approaches such as GMsFEM and NLMC
\cite{chung2018constraint,chung2018constraintmixed, NLMC}. 
We also note some works on numerical homogenization of
time fractional PDEs \cite{hu2020homogenization,brown2018numerical}.
Special constructions are needed for computing multiscale
basis functions.  These spaces satisfy the conditions
needed for partial explicit methods.


Splitting approaches are used for many applications,
\cite{marchuk1990splitting,VabishchevichAdditive}.
For example, they are used for splitting physics.
In our recent works, we have used splitting approaches 
in a design of partial explicit methods
for parabolic and wave equations
\cite{efendiev2021temporal,efendiev2021splitting}.
In the current paper, we extend these ideas to time
fractional diffusion equations, which require significant modifications.

We present several numerical results, where we compare our proposed
approaches with the approaches where all degrees of freedom
are treated implicitly. We consider high contrast permeability
fields. Our numerical results show that the proposed methods provide
similar accuracy compared to the methods where all degrees of freedom
are treated implicitly. In conclusion, we would like to highlight
some novelties of the proposed methods.

\begin{itemize}

\item The proposed methods provide a venue for performing 
partial explicit time stepping for time fractional PDEs, where
time step constraints can be severe.

\item The time step constraint depends on the coarse mesh size and
is independent of the contrast.

\item Numerical results confirm our theoretical findings.

\end{itemize}

The paper is organized as follows. In the next section, we present some
preliminaries. Section 3 is devoted to stability conditions and their derivations. In Section 4, we briefly present the construction of spaces for partially explicit method. We present the numerical results in Section \ref{sec:numerical} and conclude the paper in Section \ref{sec:conclusions}.

	\section{Preliminaries}

	Let $\Omega$ be a bounded domain in $\mathbb R^d\,(d=1,2,3)$ with
a sufficiently smooth boundary $\partial \Omega$. We consider a partial differential equation with the fractional derivative in time $t$, satisfying:
\begin{equation}
		\label{eq:model}
		\left\{
		\begin{aligned}
			\partial_t^\alpha u(x,t) & =\nabla\cdot \left(\kappa({x})\nabla u(x,t) \right) +f(x,t) && \text{in } \Omega,
			\quad t\in \left( 0,T\right] \\
			u& =0 && \text{on }\partial \Omega, \;t\in \left( 0,T\right] \\
			u(0) & =u_0(x)   && \text{in } \Omega.
		\end{aligned}
		\right .
	\end{equation}
Here, $0< \alpha <1$ is a given fixed parameter. $\kappa(x)$ is a high-contrast  multiscale field.
The initial function $u_0$ is a given term and $T>0$ is a fixed value. The source term $f$ satisfies $t^{1-\alpha}f\in L^1(0,T;L^2(\Omega))$.

In the model problem \eqref{eq:model}, $\partial_t^\alpha w$ refers to the
left-sided Caputo fractional derivative of order  $\alpha$ 
of the function $w(t)$, defined by (see, e.g. \cite[p. 91, (2.4.1)]{KilbasSrivastavaTrujillo:2006}
or \cite[p. 78]{Podlubnybook})
\begin{equation*}
  \partial^\alpha_t w (t) = \frac{1}{\Gamma(1-\alpha)}\int_0^t \frac{1}{(t-s)^\alpha}w'(s)ds.
\end{equation*}

The fractional diffusion equations were introduced in physics aiming at describing diffusions in media with fractal geometry \cite{nigmatullin1986realization}. They are enormously applied to many fields, for example in engineering, physics, biology and finance. Their practical applications include electron transport in Xerox photocopier, visco-elastic materials, and protein transport in cell membrane \cite{scher1975anomalous,giona1992fractional,kou2008stochastic}.

To discretize the model problem \eqref{eq:model}, we first decompose the time domain $[0,T]$ into $N$ time subdomains $(T_n,T_{n+1})$, $n=0,1,\cdots, N-1$, with $0=T_0<T_1<\cdots<T_{N-1}<T_{N}=T$ and $\Delta T_n:=T_{n+1}-T_{n}$. For simplicity, we assume $\Delta T_{n}=\Delta T$ for any $n=0,1,\cdots, N-1$. Let $\mathcal{T}_{H}$ be a decomposition of the spatial domain $\Omega$ into non-overlapping shape-regular rectangular elements with maximal mesh size $H$. Let $\mathcal{T}^{h}$ be a refinement of $\mathcal{T}^{H}$ with $h\ll H$. One could choose a finite element space $V_h\subset H_{0}^1(\Omega)$ with $h$ being an extremely small number and utilize full discretization to solve \eqref{eq:model}. For the sake of saving computational cost, we would like to construct a finite dimensional space $V_{H}\subset H_{0}^{1}(\Omega)$ based on $\mathcal{T}_{H}$. 

Utilizing finite difference approximation to discretize the time-fractional derivative \cite{chuanjuxu2007}, one obtains the following approximation
\begin{equation}
    \partial_t^\alpha u(x,T_k)\approx \frac{1}{\Gamma(2-\alpha)} \sum_{j=0}^{k} \frac{u(x,T_{k+1-j}-u(x,T_{k-j}))}{\Delta T^{\alpha}} [(j+1)^{1-\alpha}-j^{1-\alpha}], \text{ for } k=1,2,\cdots, N.
    \label{approx:time_fractional}
\end{equation}
Using the implicit Euler scheme and the approximation \eqref{approx:time_fractional}, one obtains the full discretized finite element method which reads as follows:
for any $k=0, 1,\cdots,N-1$, find $u^{k+1}\in V_{H}$ such that
\begin{equation}
(u^{k+1},v)+\alpha_0 a( u^{k+1}, v)=(1-b_1)(u^k,v)+\sum_{j=1}^{k-1}(u^{k-j},v)+b_k(u^0,v)+\alpha_0(f^{k+1},v), ~\forall v\in V_{H}.
\label{implicit_weak_form}
\end{equation}
Here, $\alpha_0=\Gamma(2-\alpha)\Delta T^{\alpha}$, $a(v,w):=(\kappa \nabla v,\nabla w)$ for any $v,w\in V_{H}$ and $b_j=(j+1)^{1-\alpha}-j^{1-\alpha}$, for $j=0,1,\cdots, k$.
We remark that $u^{k}$ is an approximation to the solution $u(T^{k})$.

Next, we clarify some notations used throughout the article.
We write $(\cdot,\cdot)$ to denote the inner product in $L^2(\Omega)$ and $\norm{\cdot}$ for the corresponding norm. 
Let $H^1(\Omega)$ be the classical Sobolev space with the norm $\norm{v}_1 := \left ( \norm{v}^2 + \norm{\nabla v}^2 \right )^{1/2}$ for any $v \in H^1(\Omega)$ 
and $H_0^1(\Omega)$ the subspace of functions having a vanishing trace. For any subset $S\subset \Omega$, we denote $V(S): = H_0^1(S)$.
We use $\norm{\cdot}_{a}$ to denote the norm induced by the $a$-norm. That is, $\norm{v}_a=\sqrt{a( v,v)}$.

\section{Three Schemes and Stabilities}
In this section, we will prove the unconditional stability for implicit Euler scheme in Subsection \ref{subsec:implicit}. In Subsection \ref{subsec:expilicit}, we establish the stability condition for the explicit Euler scheme. Finally, in Subsection \ref{subsec:partial_explicit}, we introduce a partial splitting algorithm and discuss the stability condition for the algorithm.

\subsection{Stability of Implicit Euler Scheme}\label{subsec:implicit}
We present in the following theorem that the implicit Euler scheme is unconditionally stable (cf. \cite{sun2020new}).

\begin{theorem}
\label{thm:implicit}
Assume the source term $f$ satisfies $t^{1-\alpha}f\in L^1(0,T;L^2(\Omega))$. Let $\Delta T$ be the time step size and define $\alpha_0:=\Gamma(2- \alpha) \Delta T^{\alpha}$. Let $u^{k}$ be the solution to \eqref{implicit_weak_form} for $k=1,2,\cdots, N$.
Then the implicit Euler Scheme \eqref{implicit_weak_form} is unconditionally stable. Moreover, we have the following stability estimate:
\[
\norm{u^{N}}^2_a\leq \norm{u^{0}}_a^2+\alpha_0\sum_{k=0}^{N-1} \norm{f^{k+1}}^2.
\]
\end{theorem}
\begin{proof}
Note that \eqref{implicit_weak_form} can be rewritten into 
\begin{equation}
\sum_{j=0}^{k}b_{k-j}\big( u^{j+1}-u^{j},v\big)+\alpha_0 (\kappa \nabla u^{k+1},\nabla v)=\alpha_0(f^{k+1},v) ~\forall v\in V_{H}.
\label{eqn:implicit_weak}
\end{equation}
Define $b_k:=b_{-k}$ for any $k<0$ and $k\in \mathbb{Z}$. Choosing $v=u^{k+1}-u^k$  in \eqref{eqn:implicit_weak} and taking a summation over $k$ from $0$ to $N-1$, one obtains the following equality.
\begin{equation}
\sum_{k=0}^{N-1} \sum_{j=0}^{k}b_{k-j}\big( u^{j+1}-u^{j},u^{k+1}-u^k\big)+\alpha_0 \sum_{k=0}^{N-1} a(u^{k+1},u^{k+1}-u^k)=\sum_{k=0}^{N-1} \alpha_0(f^{k+1},u^{k+1}-u^k).
\label{eqn:weak_sum}
\end{equation}
Notice that the first term of \eqref{eqn:weak_sum} can be written as
\[
\sum_{k=0}^{N-1} \sum_{j=0}^{k}b_{k-j}\big( u^{j+1}-u^{j},u^{k+1}-u^k\big)=\frac{b_0}{2}\sum_{k=0}^{N-1}\norm{ u^{k+1}-u^k}^2+\frac{1}{2} \sum_{k=0}^{N-1}\sum_{j=0}^{N-1}b_{|k-j|}( u^{j+1}-u^{j},u^{k+1}-u^k).
\]
Recall that 
\[
b_k=(k+1)^{1-\alpha}-k^{1-\alpha}=\frac{1}{1-\alpha} \int_{0}^{1}(k+s)^{-\alpha} ds.
\]
Notice that $\phi (k)=(k+s)^{-\alpha}$ is a complete monotonic function. It follows from Hausdorff-Bernstein-Widder Theorem that 
$(k+s)^{-\alpha} =\int_{0}^{\infty} e^{-k\tau} d g_{s} (\tau)$ for some cumulative distribution function $g_s$. We then have 
\begin{align*}
 &\sum_{k=0}^{N-1}\sum_{j=0}^{N-1}b_{|k-j|}( u^{j+1}-u^{j},u^{k+1}-u^k)   \\
 =&\frac{1}{1-\alpha} \int_{0}^{1}\int_{0}^{\infty}\sum_{k=0}^{N-1}\sum_{j=0}^{N-1} e^{-|k-j|\tau} ( u^{j+1}-u^{j})(u^{k+1}-u^k) d g_s(\tau) ds.
\end{align*}
Create a matrix $M$ of size $N\times N$ with entry $M_{k,j}=e^{-|k-j|t}$ for $0\leq k,j\leq N-1$.
Then the matrix $M$ is positive definite for $t>0$. Therefore, 
\[
\sum_{k=0}^{N-1}\sum_{j=0}^{N-1}b_{|k-j|}( u^{j+1}-u^{j}),u^{k+1}-u^k)>0.
\]
Therefore, we have 
\[
\sum_{k=0}^{N-1} \sum_{j=0}^{k}b_{k-j}\big( u^{j+1}-u^{j},u^{k+1}-u^k\big)+\alpha_0 \sum_{k=0}^{N-1} a(u^{k+1},u^{k+1}-u^k)\geq \frac{1}{2} \sum_{k=0}^{N-1} \norm{u^{k+1}-u^{k}}^2+\alpha_0 \sum_{k=0}^{N-1} a(u^{k+1},u^{k+1}-u^{k}).
\]
Notice that $a(u^{k+1},u^{k+1}-u^{k})=\frac{1}{2}(\norm{u^{k+1}}^2_{a}-\norm{u^{k}}^2_{a}+\norm{u^{k+1}-u^k}^2_{a}) $.
Furthermore, utilizing Cauchy-Schwartz Inequality, one obtains  the following estimate.
\begin{align}
\alpha_0 \sum_{k=0}^{N-1} (f^{k+1},u^{k+1}-u^{k})&\leq \sum_{k=0}^{N-1} \alpha_0\norm{f^{k+1}}\norm{ u^{k+1}-u^{k}}\\
    &\leq \frac{\alpha_0^2}{2}\sum_{k=0}^{N-1}    \norm{f^{k+1}}^2+\frac{1}{2}\sum_{k=0}^{N-1}  \norm{ u^{k+1}-u^{k}}^2
\end{align}
Therefore we have 
$\frac{1}{2} \sum_{k=0}^{N-1} \norm{u^{k+1}-u^{k}}^2+\sum_{k=0}^{N-1} \frac{\alpha_0}{2} (\norm{u^{k+1}}^2_{a}-\norm{u^{k}}^2_{a}+\norm{u^{k+1}-u^k}^2_{a}) \leq  \frac{\alpha_0^2}{2}\sum_{k=0}^{N-1}   \norm{f^{k+1}}^2+\frac{1}{2}\sum_{k=0}^{N-1}  \norm{ u^{k+1}-u^{k}}^2.$ It follows from $f$ satisfies $t^{1-\alpha}f\in L^1(0,T;L^2(\Omega))$ that $\alpha_0\sum_{k=0}^{N-1} \norm{f^{k+1}}^2$ is bounded.
So we have 
\begin{equation}
    \norm{u^{N}}^2_a\leq \norm{u^{0}}_a^2+\alpha_0\sum_{k=0}^{N-1} \norm{f^{k+1}}^2.
\label{estimate:implicit}
\end{equation}
\end{proof}

\subsection{Stability of Explicit Euler Scheme}\label{subsec:expilicit}
With explicit Euler scheme and the approximation \eqref{approx:time_fractional}, the full discretized finite element method reads as follows:
for any $k=0, 1,\cdots,N-1$, find $u^{k+1}\in V_{H}$ such that
\begin{equation}
(u^{k+1},v)+\alpha_0 a( u^{k}, v)=(1-b_1)(u^k,v)+\sum_{j=1}^{k-1}(u^{k-j},v)+b_k(u^0,v)+\alpha_0(f^{k+1},v), ~\forall v\in V_{H}.
\label{scheme:explicit}
\end{equation}

By a similar argument and $a(u^k,u^{k+1}-u^k)=\frac{1}{2}(\norm{u^{k+1}}_{a}^{2}-\norm{u^{k}}_{a}^{2}-\norm{u^{k+1}-u^{k}}_{a}^{2})$, we present the stability condition of explicit Euler scheme in the following theorem. For the brevity of the paper, we omit the proof here.
\begin{theorem}
Assume the source term $f$ satisfies $t^{1-\alpha}f\in L^1(0,T;L^2(\Omega))$. Let $\Delta T$ be the time step size and define $\alpha_0:=\Gamma(2- \alpha) \Delta T^{\alpha}$. Let $u^{k}$ be the solution to \eqref{scheme:explicit} for $k=1,2,\cdots, N$.
Then the explicit scheme \eqref{scheme:explicit} is stable
if $\norm{v}^2\geq 2\alpha_0 \norm{v}_{a}^{2}$ for any $v\in V_{H}$. Moreover, we
have the following estimate
\begin{equation}
\norm{u^{N}}_a^2\leq \norm{u^{0}}_a^2+\alpha_0\sum_{k=0}^{N-1}\norm{f^{k+1}}^2.
\label{stability_explicit}
\end{equation}
\end{theorem}

\subsection{Stability of Partially Explicit Scheme}\label{subsec:partial_explicit}
In this subsection, we shall introduce a partially explicit temporal splitting scheme. For this purpose, we decompose the solution space $V_{H}$ into two subspaces $V_{H,1}$ and $V_{H,2}$. That is, $V_{H}=V_{H,1}+V_{H,2}$. Then
 the solution $u^{k+1}$ can be written as $u^{k+1}=u^{k+1}_{1}+u^{k+1}_{2}$ with $u^{k+1}_1\in V_{H,1}$ and $u^{k+1}_{2}\in V_{H,2}$. 
 The full discretized finite element method with partially explicit temporal splitting scheme reads as follows:
for any $k=0, 1,\cdots,N-1$, find $u^{k+1}_{1}\in V_{H,1}$ and $u^{k+1}_{2}\in V_{H,2}$ such that
\begin{equation}
 (u^{k+1},v)+\alpha_0 a(u^{k+1}_{1}+ u^{k}_{2}, v)=(1-b_1)(u^k,v)+\sum_{j=1}^{k-1}(u^{k-j},v)+b_k(u^0,v)+\alpha_0(f^{k+1},v), ~\forall v\in V_{H}.
\label{scheme:splitting}   
\end{equation}
We study the stability of the partially explicit scheme in the following theorem.
\begin{theorem}
Assume the source term $f$ satisfies $t^{1-\alpha}f\in L^1(0,T;L^2(\Omega))$. Let $\Delta T$ be the time step size and define $\alpha_0:=\Gamma(2- \alpha) \Delta T^{\alpha}$. Let $u^{k}$ be the solution to \eqref{scheme:splitting} for $k=1,2,\cdots, N$.
Then the partially explicit scheme \eqref{scheme:splitting} is stable if 
\begin{equation}
\norm{u_{1}+u_{2}}^{2}\geq 2(1-\gamma^{2})\norm{u_2}^2    \text{ for any } u_{1}\in V_{H,1}, u_{2}\in V_{H,2}
\label{Splitting_cond1}
\end{equation}
and
\begin{equation}
  \norm{v}\geq \alpha_0 (1-\gamma^2)^{-1}\norm{v}_{a}^{2} \text{ for any } v\in V_{H,2},
  \label{Splitting_cond2}
\end{equation}
 for some constant $\gamma\in[0,1)$. Furthermore, 
we  have the following estimate
 \begin{equation}
     \norm{u^{N}}_a^2\leq \norm{u^{0}}_a^2+\sum_{k=0}^{N-1}\alpha_0\norm{f^{k+1}}^2.
 \end{equation}
\end{theorem}
\begin{proof}
Notice that
\begin{align*}
    a(u^{k+1}_{1}-u^{k}_{2},u^{k+1}_{}-u^{k}_{})&=a(u^{k+1},u^{k+1}_{}-u^{k}_{})+a(u^{k}_{2}-u^{k+1}_{2},u^{k+1}_{}-u^{k}_{})\\
    &\geq \frac{1}{2} (\norm{u^{k+1}}_{a}^{2}-\norm{u^{k}}_{a}^{2}-\norm{u^{k+1}-u^{k}}_{a}^{2})-\frac{1}{2} (\norm{u^{k+1}_{2}-u^{k}_{2}}_{a}^{2}+\norm{u^{k+1}_{}-u^{k}_{}}_{a}^{2})\\
    &=\frac{1}{2} (\norm{u^{k+1}}_{a}^{2}-\norm{u^{k}}_{a}^{2}-\norm{u^{k+1}_{2}-u^{k}_{2}}_{a}^{2}).
\end{align*}

 Applying a similar argument adopted in the proof of Theorem \ref{thm:implicit} and combining \eqref{Splitting_cond1} and \eqref{Splitting_cond2}, we  have 
 \begin{equation}
     \norm{u^{N}}_a^2\leq \norm{u^{0}}_a^2+\sum_{k=0}^{N-1}\alpha_0\norm{f^{k+1}}^2.
 \end{equation}
 
\end{proof}

\section{$V_{H,1}$ and $V_{H,2}$ Constructions}

\begin{figure}[H]
		\centering
		\includegraphics[trim={2.5cm 1.5cm 2.0cm 0.2cm},clip,width=0.35 \textwidth]{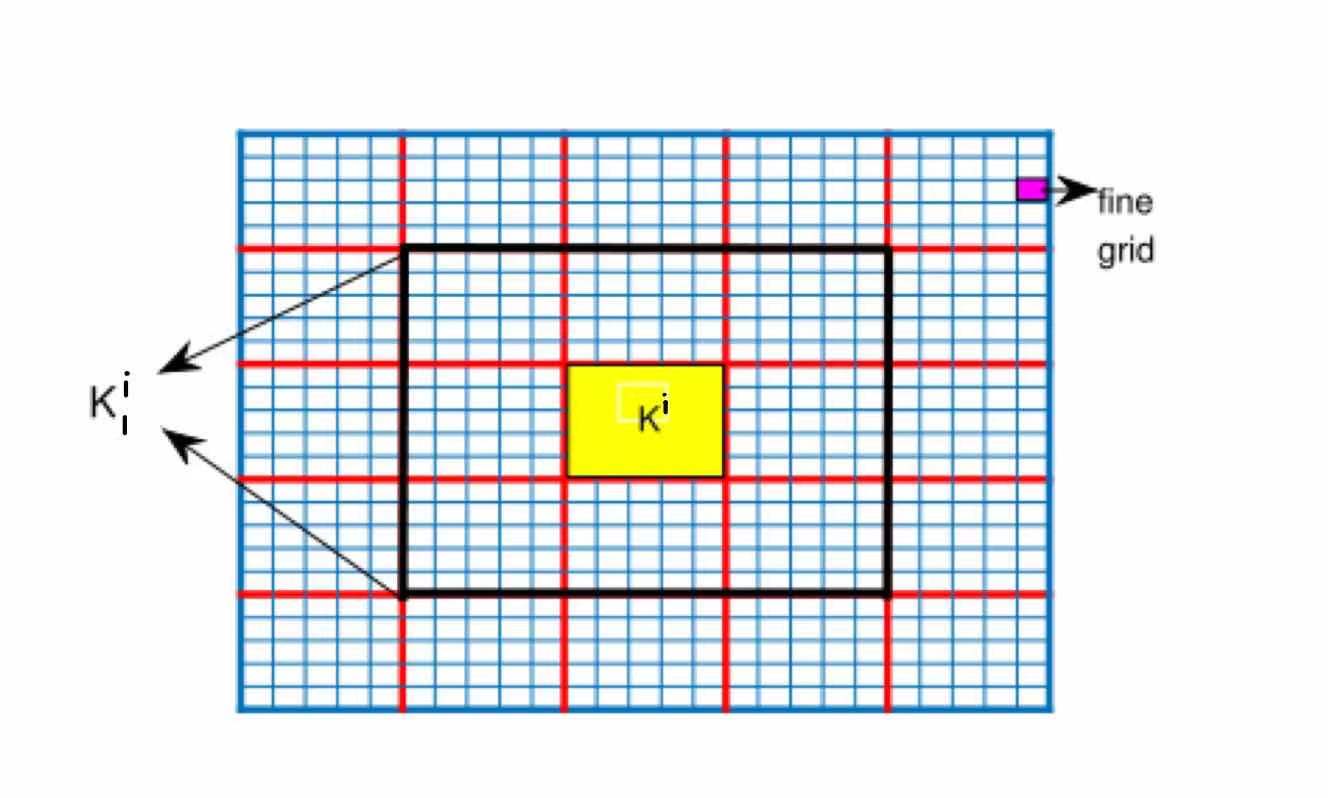}
		\caption{Illustration of the oversampling domain.}
		\label{fig:oversamp}
	\end{figure}

In this section, we shall introduce a possible way to construct spaces $V_{H,1}$ and $V_{H,2}$ such that the partially explicit scheme is stable.
Our construction follows our previous work \cite{chung2021contrastindependent}.
 We will first recap the constrained  energy minimization (CEM) finite element method and show that the
CEM type finite element space is a good choice
of $V_{H,1}$ since the CEM multiscale basis functions are constructed in a way that
they are almost orthogonal to the space $\tilde{V}$ which will be defined in Section \ref{sec:cem}. We will present two possible ways of constructing the subspace $V_{H,2}$. Before that, we introduce a concept called oversampling domain, which will be used later. For each coarse element $K_i$, we define the oversampled region $K_{i,k_i} \subseteq \Omega$ by enlarging $K_i$ by $k_i \in \mathbb{N}$ layer(s), i.e., 
$$ K_{i,0} := K_i, \quad K_{i,k_i} := \bigcup \{ K \in \mathcal{T}^H : K \cap K_{i,k_i -1} \neq \emptyset \} \quad \text{for } k_i = 1, 2, \cdots.$$
We call $k_i$ a parameter of oversampling related to the coarse element $K_i$. See Figure \ref{fig:oversamp} for an illustration of $K_{i,1}$. For simplicity, we denote $K_i^+$ a generic oversampling region related to the coarse element $K_i$ with a specific oversampling parameter $k_i$. 

\subsection{CEM method}
\label{sec:cem}

In this section, we will discuss the CEM method for solving the problem
\eqref{implicit_weak_form}. 
In particular, we will focus on constructing the finite element space by solving
a constrained energy minimization problem. Let $\mathcal{T}_{H}$ be a decomposition of the spatial domain $\Omega$ into non-overlapping shape-regular rectangular elements with maximal mesh size $H$. We shall construct a finite dimensional space $V_{H}\subset H_{0}^{1}(\Omega)$. For each element $K_{i}\in\mathcal{T}_{H}$, one obtains a set of auxiliary basis functions $\{\psi_{j}^{(i)}\}_{j=1}^{L_{i}}\subset V(K_{j})$ by solving eigenvalue problems: finding $(\lambda_{j}^{(i)},\psi_{j}^{(i)}) \in \mathbb{R}\times V(K_{j}) $ such that
\begin{equation*}
\int_{K_{i}} \kappa \nabla \psi_{j}^{i} \cdot \nabla v= \lambda_{j}^{i} s_i(\psi_{j}^{i},v), ~~\forall v\in V(K_{i}),
\end{equation*}
where 
$s_{i}(u,v)=\int_{K_{i}}\tilde{\kappa}uv$
and  $\tilde{\kappa}=\kappa\sum_{i=1}^{N_c}|\nabla\chi_{i}|^{2}$
with $\{\chi_{i}\}_{i=1}^{N_c}$ being the standard multiscale finite element basis functions. Select the first $L_i$ eigenfunctions corresponding to the first small eigenvalues. $\{\psi_{j}^{(i)}\}_{j=1}^{L_{i}}$ is the set of the auxiliary basis functions and we denote $V_{aux}^{(i)}:=\text{span}\{\psi_{j}^{(i)}:\;1\leq j\leq L_{i}$ as the local auxiliary space.
We then can define a local projection operator $\Pi_{K_{i}}:L^{2}(K_{i})\mapsto V_{aux}^{(i)}\subset L^{2}(K_{i})$
such that 
\[
s_{i}(\Pi_{i}u,v)=s_{i}(u,v)\;\forall v\in V_{aux}^{(i)}.
\]

Next we define a global  projection operator by $\Pi:L^{2}(\Omega)\mapsto V_{aux}\subset L^{2}(\Omega)$ such that
\[
s(\Pi u,v)=s(u,v)\;\forall v\in V_{aux}:=\bigcup_{i=1}^{N_{e}}V_{aux}^{(i)},
\]
where $s(u,v):=\sum_{i=1}^{N_{e}}s_{i}(u|_{K_{i}},v|_{K_{i}})$ and $N_e$ is the number of coarse elements.
For each auxiliary basis functions $\psi_{j}^{(i)}$, we can define
a local basis function $\phi_{j}^{(i)}\in V(K_{i}^{+})$
such that 
\begin{align*}
a(\phi_{j}^{(i)},v)+s(\mu_{j}^{(i)},v) & =0& \forall v\in V(K_{i}^{+})\\
s(\phi_{j}^{(i)},\nu) & =s(\psi_{j}^{(i)},\nu)& \forall\nu\in V_{aux}(K_{i}^{+}).
\end{align*}
We then define the space $V_{cem}$ as 
\begin{align*}
V_{cem} & :=\text{span}\{\phi_{j}^{(i)}:\;1\leq i\leq N_{e},1\leq j\leq L_{i}\}.
\end{align*}
The CEM solution $u^{k+1}_{cem}$ for $k=1,2,\cdots,N-1$ is given by
\begin{equation}
(u_{cem}^{k+1},v)+\alpha_0 a( u_{cem}^{k+1}, v)=(1-b_1)(u_{cem}^k,v)+\sum_{j=1}^{k-1}(u_{cem}^{k-j},v)+b_k(u_{cem}^0,v)+\alpha_0(f^{k+1},v), ~\forall v\in V_{cem}.
\label{scheme:CEM}
\end{equation}
Now we construct global basis functions $\phi_{j,glo}^{(i)}$. For each auxiliary basis functions $\psi_{j}^{(i)}$, we find $\phi_{j,glo}^{(i)}\in V$ 
such that 
\begin{align*}
a(\phi_{j,glo}^{(i)},v)+s(\mu_{j}^{(i)},v) & =0& \forall v\in V\\
s(\phi_{j,glo}^{(i)},\nu) & =s(\psi_{j,glo}^{(i)},\nu)& \forall\nu\in V_{aux}.
\end{align*}
We remark here that the local multiscale basis $\phi_{j}^{(i)}$ is an approximation of the global basis function $\phi_{j,glo}^{(i)}$. Denote $V_{glo}:=\text{span}\{\phi_{j,glo}^{(i)}:\;1\leq i\leq N_{e},1\leq j\leq L_{i}\}.$
It can be proved the $V_{glo}$ is $a-$orthogonal to a space $\tilde{V}:=\{v\in V:\;\Pi(v)=0\}$.
We also know that $V_{cem}$ is closed in $V_{glo}$ and therefore
it is almost orthogonal to $\tilde{V}$. Thus, we can choose $V_{cem}$
to be $V_{H,1}$ and it remains to construct a space $V_{H,2}$ in $\tilde{V}$.

\subsection{Construction of $V_{H,2}$}
\label{subsection:second_choice}
In this subsection, we will discuss a choice for the space $V_{H,2} \subset \tilde{V}$. This choice of $V_{H,2}$ is based on the CEM type multiscale finite
element space. For
each coarse element $K_{i}$, we will solve an eigenvalue problem to obtain the auxiliary
basis. Find eigenpairs $(\xi_{j}^{(i)},\gamma_{j}^{(i)})\in(V(K_{i})\cap\tilde{V})\times\mathbb{R}$ such that
\begin{align}
\label{eq:spectralCEM2}
\int_{K_{i}}\kappa\nabla\xi_{j}^{(i)}\cdot\nabla v & =\gamma_{j}^{(i)}\int_{K_{i}}\xi_{j}^{(i)}v, \;\ \forall v\in V(K_{i})\cap\tilde{V}.
\end{align}
For each $K_i$, we order the eigenvalues in an increasing order and  choose the first smallest few $J_i$ eigenfunctions corresponding to the  $J_i$ eigenvalues. The auxiliary space $V_{aux,2}$ is formed by the span of these $J_i$ eigenfunctions. We use the notation $V_{aux,1}$ to denote the space $V_{aux}$ defined in Section \ref{sec:cem}.
For each basis function $\xi_j^{(i)} \in V_{aux,2}$, we will define a basis function $\zeta_{j}^{(i)}$: find
$(\zeta_{j}^{(i)}, \mu_{j}^{(i)},\mu_{j}^{(i),2} ) \in  V(K_i^+)\times  \in V_{aux,1}\times   V_{aux,2}$ such that 
\begin{align}
a(\zeta_{j}^{(i)},v)+s(\mu_{j}^{(i),1},v)+ ( \mu_{j}^{(i),2},v) & =0, \;\forall v\in V(K_i^+), \label{eq:v2a} \\
s(\zeta_{j}^{(i)},\nu) & =0, \;\forall\nu\in V_{aux,1}, \label{eq:v2b} \\
(\zeta_{j}^{(i)},\nu) & =( \xi_{j}^{(i)},\nu), \;\forall\nu\in V_{aux,2}. \label{eq:v2c}
\end{align}
We define $$V_{H,2}:=\text{span}\{\zeta_{j}^{(i)}| \; \forall K_i, \; \forall 1 \leq j\leq J_i\}.$$

\section{Numerical Experiment} \label{sec:numerical}
In this section, we shall present numerical results to demonstrate the performance of our proposed partially explicit scheme to solve the time-fractional diffusion equations. We  consider the time-fractional diffusion equation \eqref{eq:model} in the unit square $\Omega :=[0,1]^2$ with the final time $T:=0.01$.  The time mesh size is chosen as $\Delta T:=2*10^{-5}$ to discretize the time domain.
Let $\mathcal{T}_H$ be a decomposition of the spatial domain $\Omega$ into non-overlapping shape-regular rectangular elements with maximal mesh size $H:=\frac{1}{10}$. Since there is no analytic solution to system \eqref{eq:model}, we need to find an approximation of the exact solutions. To this end, the coarse rectangular elements are further partitioned into a collection of connected fine rectangular elements $\mathcal{T}_h$ using fine mesh size $h:=\frac{1}{100}$. Similarly, we define $V_h$ to be a conforming piecewise affine finite element associated with $\mathcal{T}_h$.\\
To ensure the fine solutions better approximation to the exact solutions, we further partition the time mesh into the fine time mesh with the mesh size $\delta T=\frac{\Delta T}{5}$.
We will use the constructed fine spatial mesh, fine time mesh and conforming finite element method to obtain the reference solutions $U_h^{k+1}$: 
for any $k=0, 1,\cdots,\frac{T}{\delta T}-1$, find $U_{h}^{k+1}\in V_{h}$ such that
\begin{equation*}
(U_{h}^{k+1},v)+\alpha_0 a( U_{h}^{k+1}, v)=(1-b_1)(U_{h}^k,v)+\sum_{j=1}^{k-1}(U_{h}^{k-j},v)+b_k(U_{h}^0,v)+\alpha_0(f^{k+1},v), ~\forall v\in V_{h}.
\end{equation*}
Notice that $U^{k}_{h}$ is an approximation of $u(x,T_{k})$ for $k=0, 1,\cdots,N-1, N$. 
In our numerical experiments, the space meshes $\mathcal{T}_H$ and $\mathcal{T}_h$ will be fixed.  

To observe the performance of the partially explicit scheme, we present three numerical solutions in addition to the fine solutions.
We use $V_{cem}$ as a solution space and seek for the CEM solutions $u^{k}_{cem}$ by \eqref{scheme:CEM}, for $k=1,2,\cdots, \frac{T}{\Delta T}$. The second numerical solutions are sought in the solution space $V_{cem}+V_{H,2}$ with implicit scheme, where $V_{H,2}$ is the space constructed in Subsection \ref{subsection:second_choice}.  The scheme reads: for any $k=0, 1,\cdots,\frac{T}{\Delta T}-1$, find $\tilde{U}^{k+1}\in V_{cem}+V_{H,2}$ such that
\begin{equation*}
(\tilde{U}^{k+1},v)+\alpha_0 a( \tilde{U}^{k+1}, v)=(1-b_1)(\tilde{U}^k,v)+\sum_{j=1}^{k-1}(\tilde{U}^{k-j},v)+b_k(\tilde{U}^0,v)+\alpha_0(f^{k+1},v), ~\forall v\in V_{cem}+V_{H,2}.
\end{equation*}
The last numerical solutions $u_{scem}^{k}$ are obtained using the partially explicit scheme \eqref{scheme:splitting}, for $k=1,2,\cdots, \frac{T}{\Delta T}$.

Our numerical experiments include testing smooth source term in Subsection \ref{subsec:numerical_smooth}
and discontinuous source term in Subsection \ref{subsection:numerical_discontinuous}.
\subsection{Numerical Experiment 1: smooth source term}\label{subsec:numerical_smooth}
\begin{figure}[H]
		\centering
		\includegraphics[trim={2.5cm 1.5cm 2.0cm 0.2cm},clip,width=0.30 \textwidth]{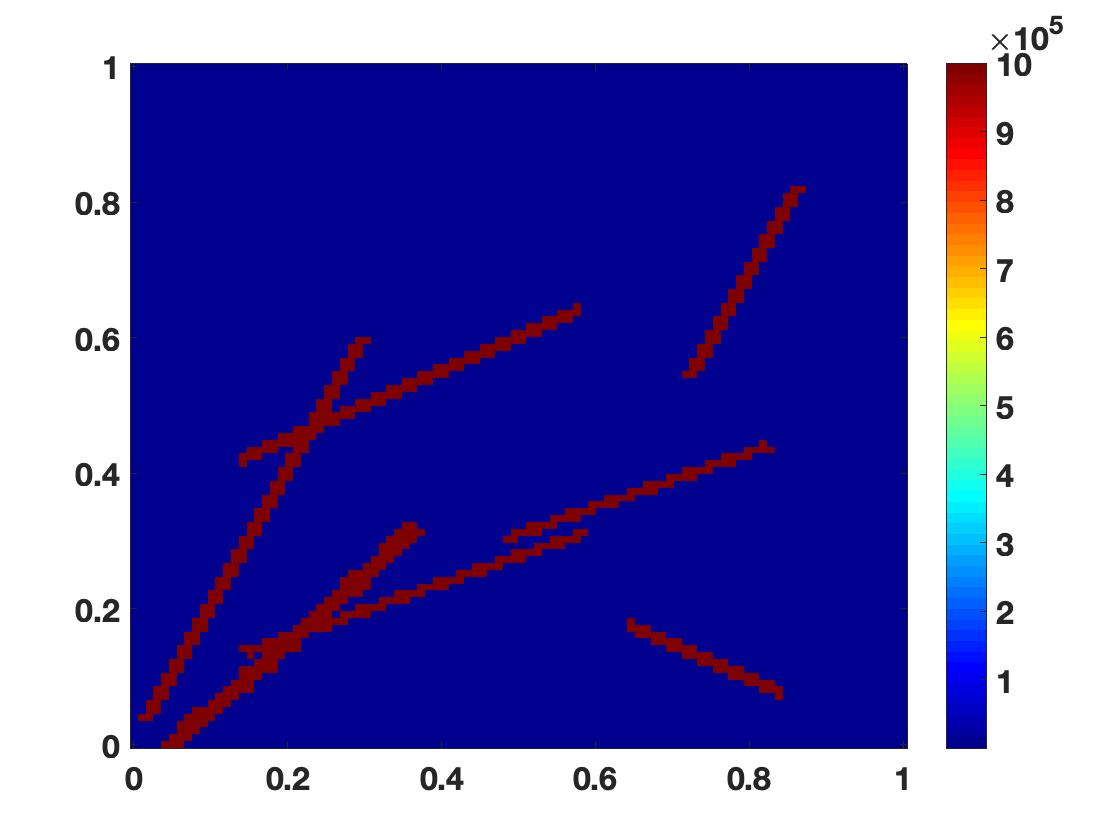}
		\includegraphics[trim={2.5cm 1.5cm 2.0cm 0.2cm},clip,width=0.30 \textwidth]{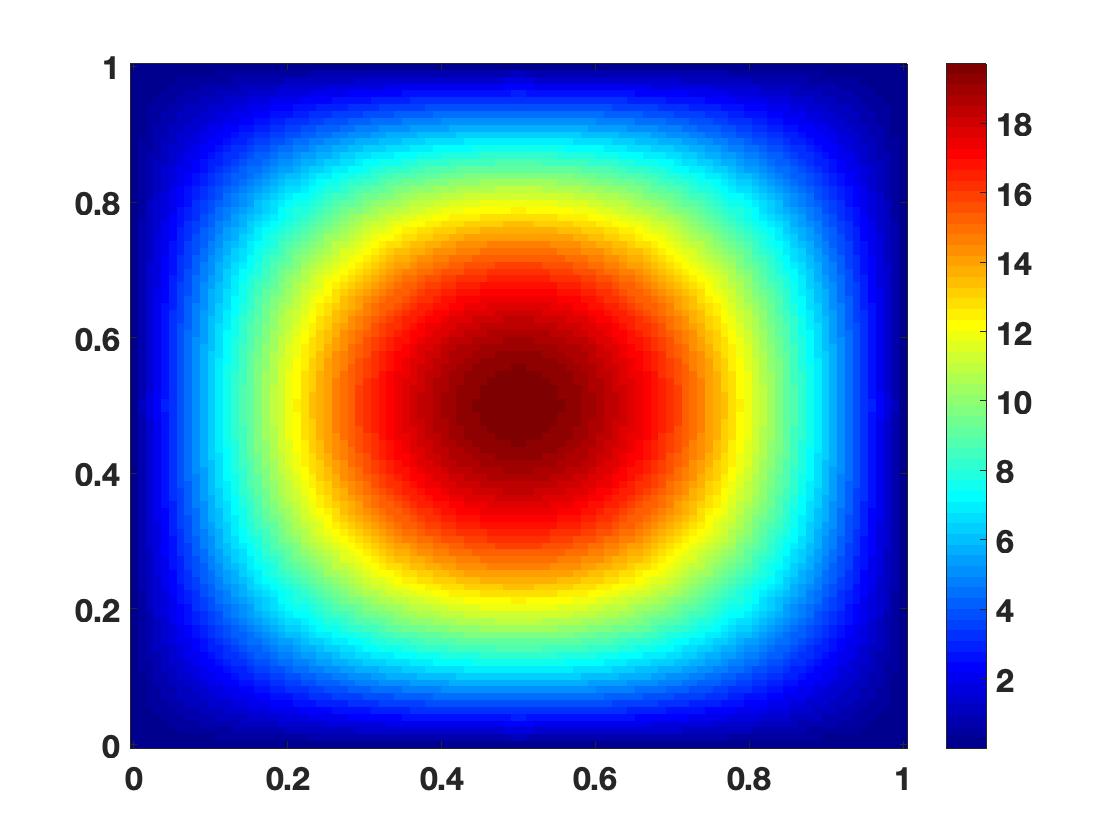}
		\caption{The heterogeneous permeability field: $\kappa$ (left) and source term $f=2\pi ^2\sin(\pi x_1)\sin(\pi x_2)$ (right).}
		\label{fig:kappa_1}
	\end{figure}

\begin{figure}[H]
		\centering
		\includegraphics[trim={2.5cm 1.5cm 1.0cm 0.2cm},clip,width=0.28 \textwidth]{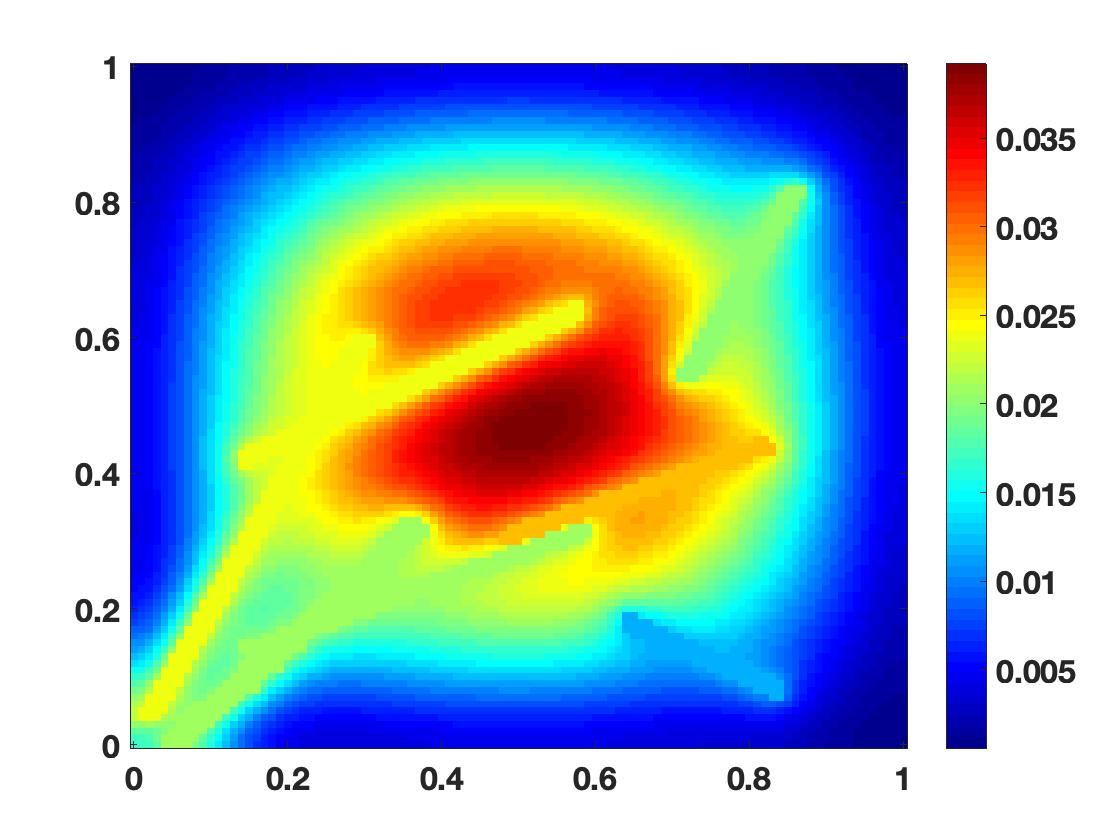}
		\includegraphics[trim={2.5cm 1.5cm 1.0cm 0.2cm},clip,width=0.28 \textwidth]{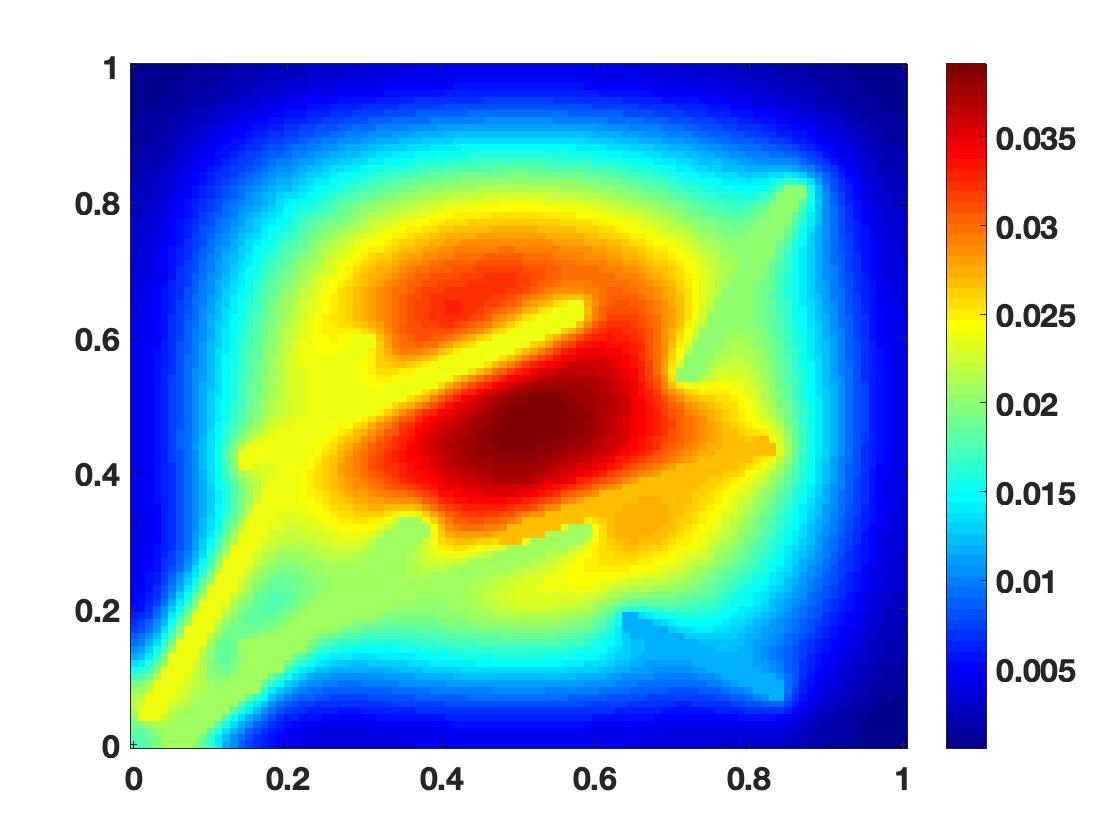}\\
		\includegraphics[trim={2.5cm 1.5cm 1.0cm 0.2cm},clip,width=0.28 \textwidth]{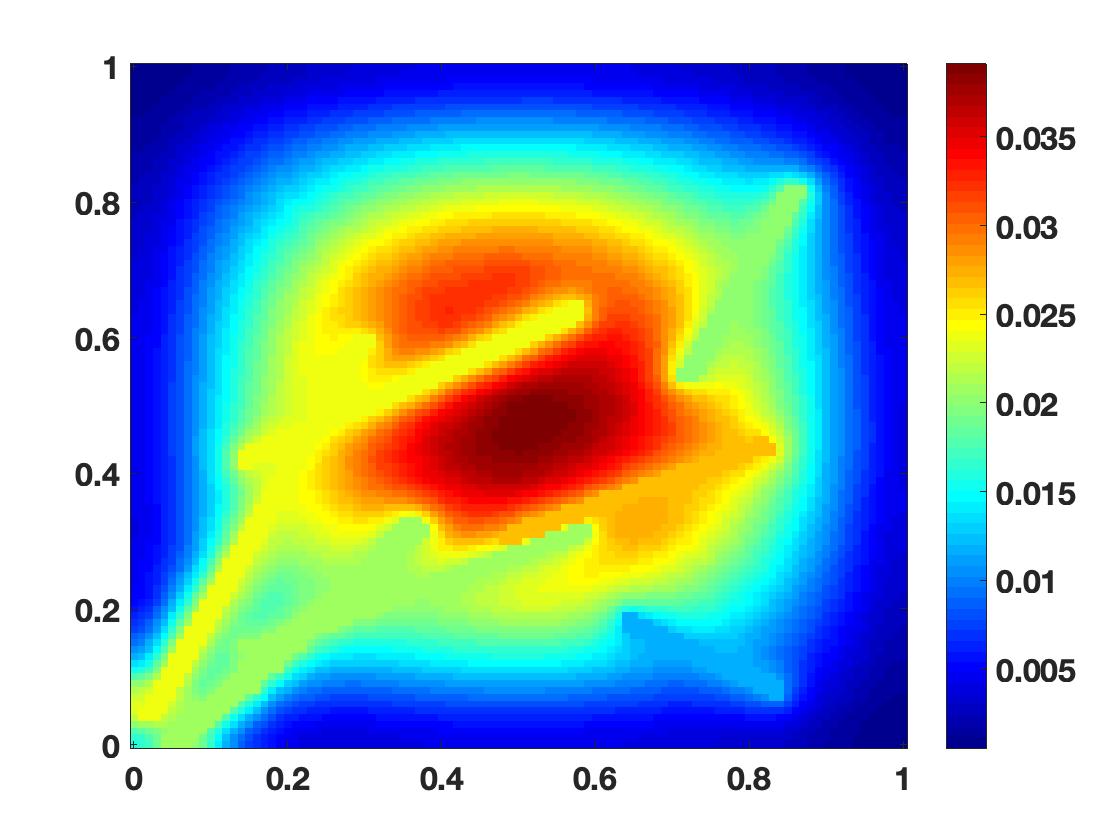}
		\includegraphics[trim={2.5cm 1.5cm 1.0cm 0.2cm},clip,width=0.28 \textwidth]{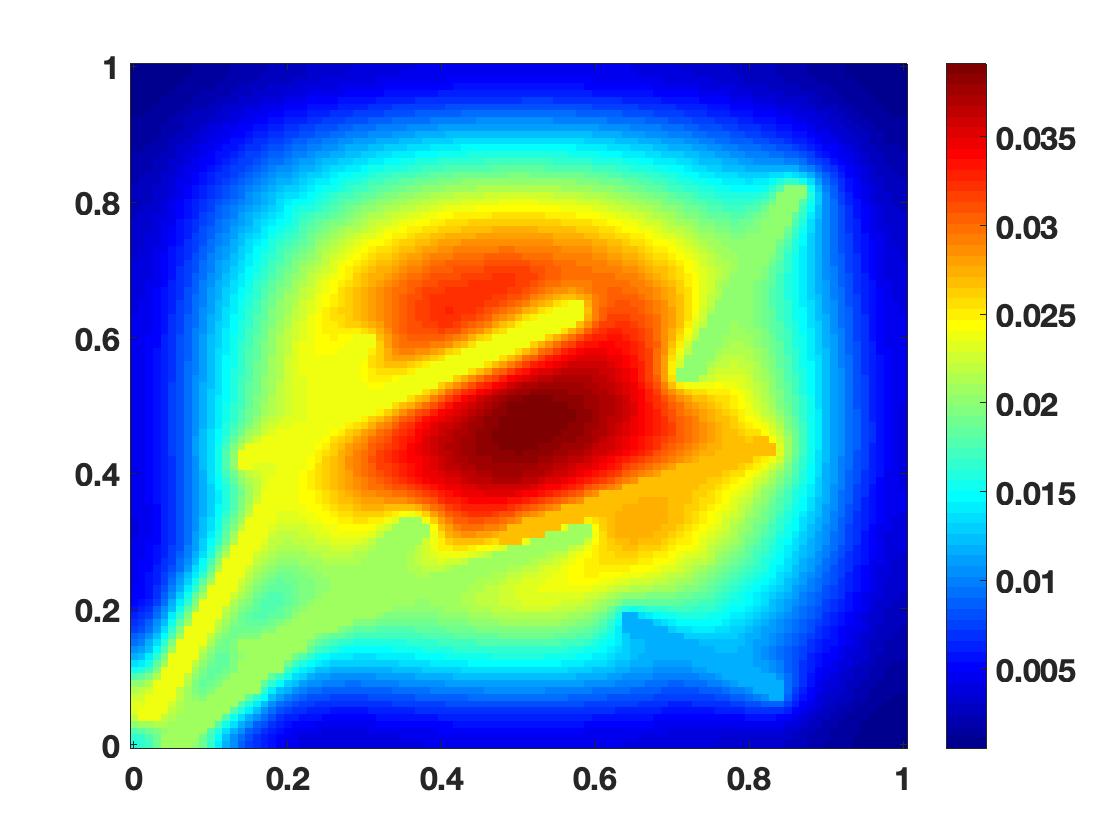}
		\caption{fine-grid solution $u_{h}^{2500}$ (top left), CEM solution $U_{cem}^{500}$ (top right), CEM solution with additional basis functions implicitly $\tilde{U}^{500}$ (bottom left), SCEM solution with additional basis functions explicitly $u_{scem}^{500}$  (bottom right).}
		\label{fig:solution_test1}
	\end{figure}

In this experiment, we choose a heterogeneous permeability coefficient $\kappa$, which has two distinct value: 1 and $10^5$. The source term is chosen to be a smooth function $f(x_1,x_2,t):=2\pi ^2\sin(\pi x_1)\sin(\pi x_2)$.
 The permeability field and the source term are plotted in Figure \ref{fig:kappa_1} for an illustration.
The fractional derivative order is chosen to be $\alpha:=0.9$. For the brevity of the paper, we will only present numerical solutions at the final time $T$. The fine-grid solution, CEM solution, CEM solution with more basis functions and SCEM solution at the final time $T$ are plotted in Figure \ref{fig:solution_test1}.
	
\begin{figure}[H]
		\centering
		\includegraphics[trim={2.5cm 0.0cm 1.0cm 0.2cm},clip,width=0.35 \textwidth]{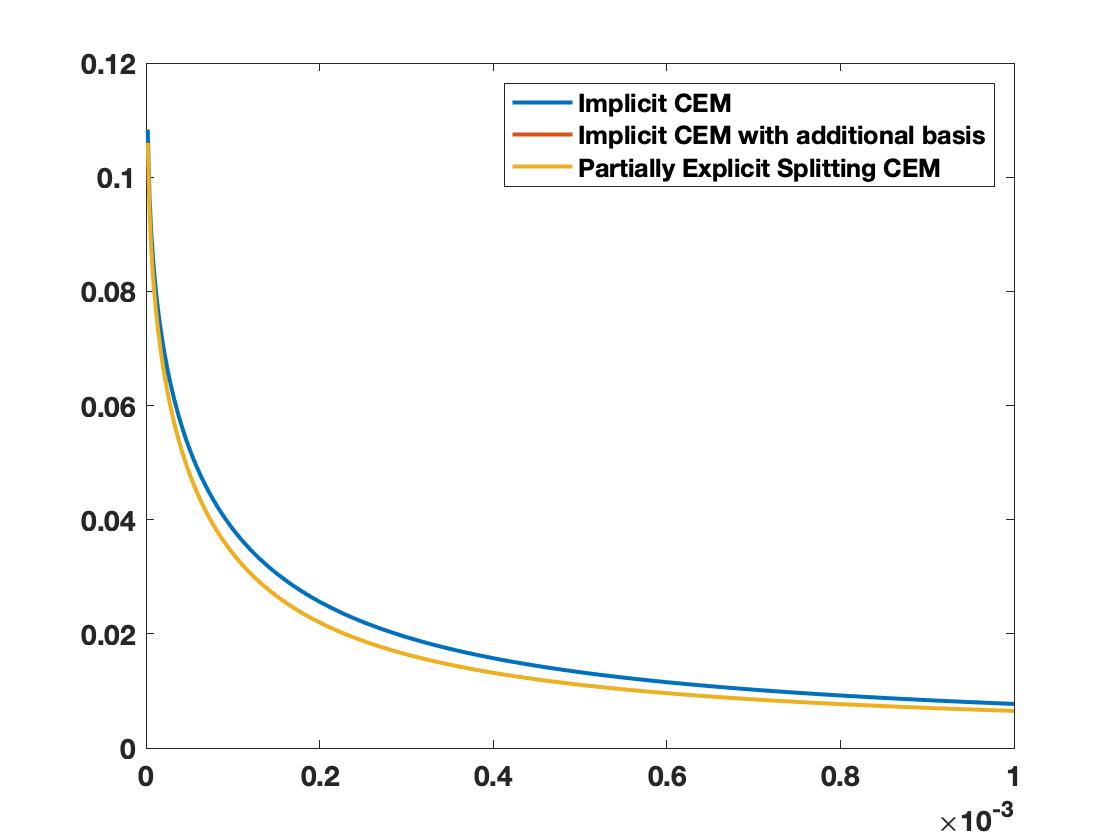}
		\includegraphics[trim={2.5cm 0.0cm 1.0cm 0.2cm},clip,width=0.35 \textwidth]{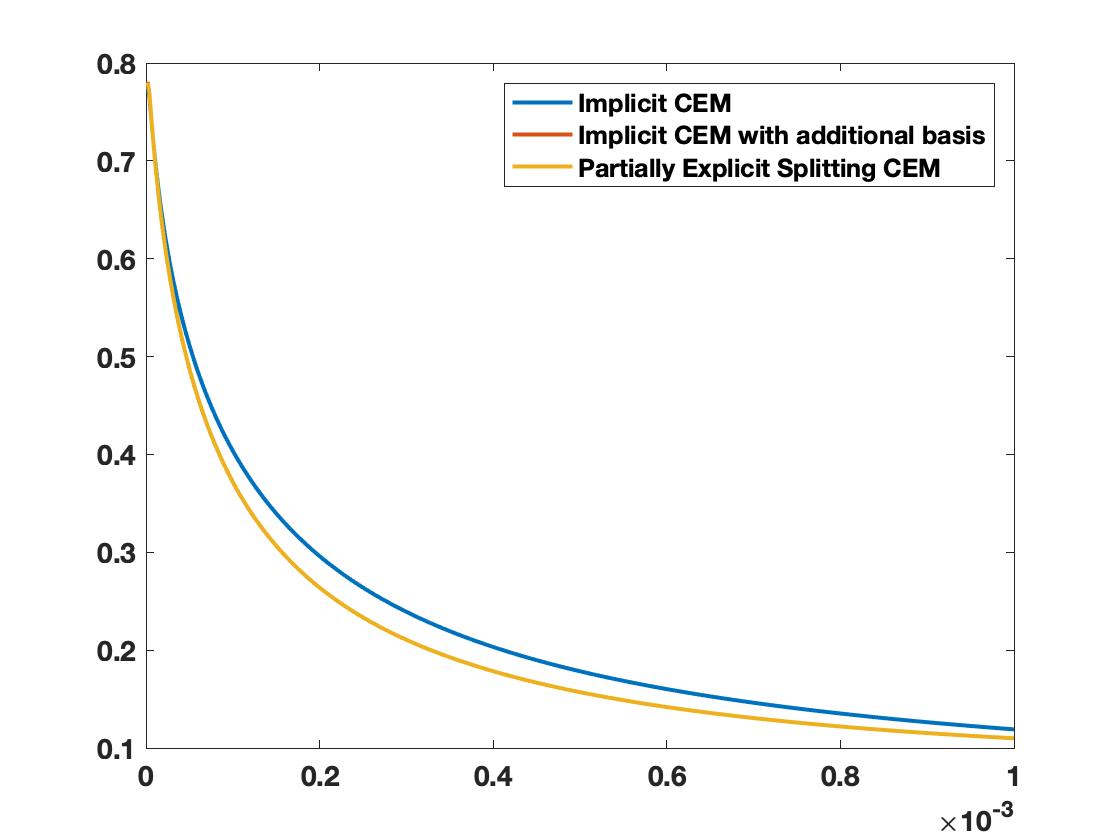}
		\caption{Relative $L^2(\Omega)$ error (left) and Relative $H^1_{\kappa}(\Omega)$ error (right).}
		\label{fig:error_test1}
	\end{figure}

The convergence history of three numerical solutions in relative $L^2(\Omega)$-norm and relative $H^1_{\kappa}(\Omega)$-norm are presented 
 in Figure \ref{fig:error_test1}. 
From Figure \ref{fig:error_test1}, one can see that when using the same number of basis functions, SCEM solutions are better approximations than CEM solutions to the reference solutions.	Moreover, numerical solutions $u_{scem}$ and $\tilde{U}$ have about the similar accuracy.
However, it is computationally cheaper to solve for $u_{scem}$ than $\tilde{U}$.

Furthermore, we test the experiment with different value of fractional derivative order $\alpha$.  It turns out that when we choose $\alpha=0.5, 0.4, 0.3$, the SCEM solutions become unstable. This observation is confirmed by the stability condition that $\norm{v}\geq \alpha_0 (1-\gamma^2)^{-1}\norm{v}_{a}^{2}$ hold true for any $v\in V_{H,2}$.
\subsection{Numerical Experiment 2: discontinuous source term}\label{subsection:numerical_discontinuous}
\begin{figure}[H]
		\centering
		\includegraphics[trim={2.5cm 1.5cm 2.0cm 0.2cm},clip,width=0.3 \textwidth]{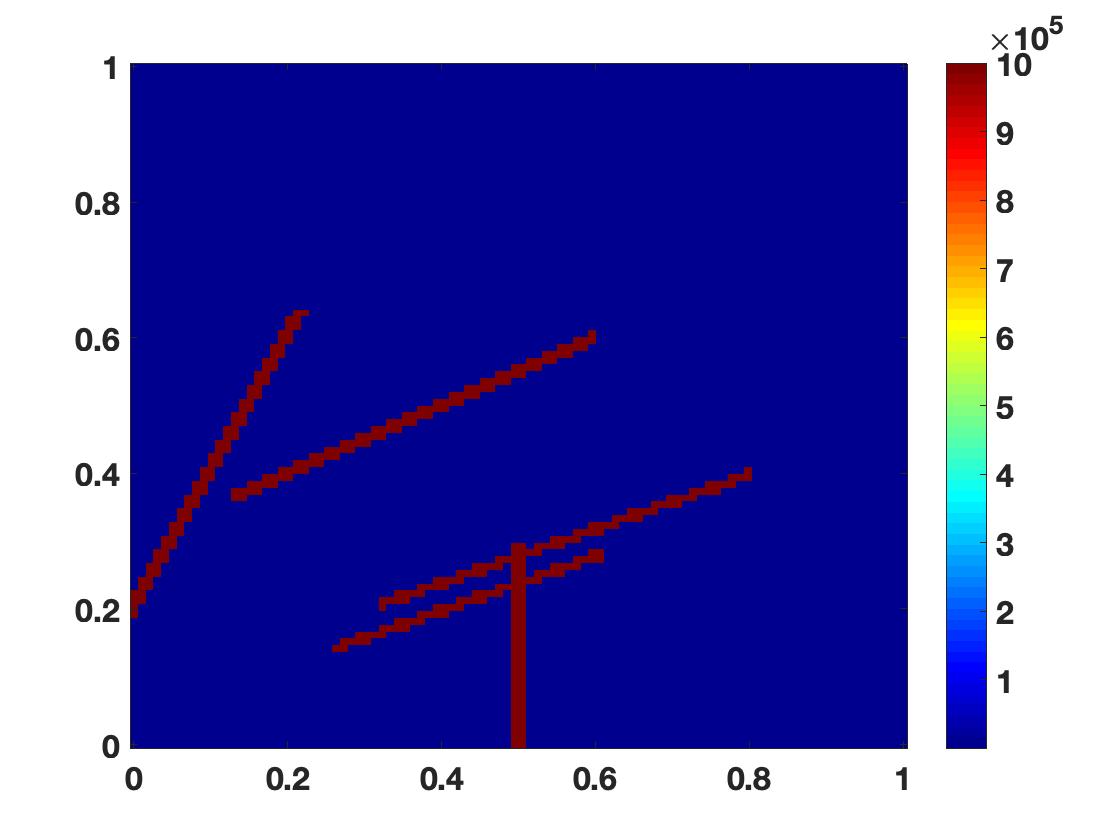}
		\includegraphics[trim={2.5cm 1.5cm 2.0cm 0.2cm},clip,width=0.3 \textwidth]{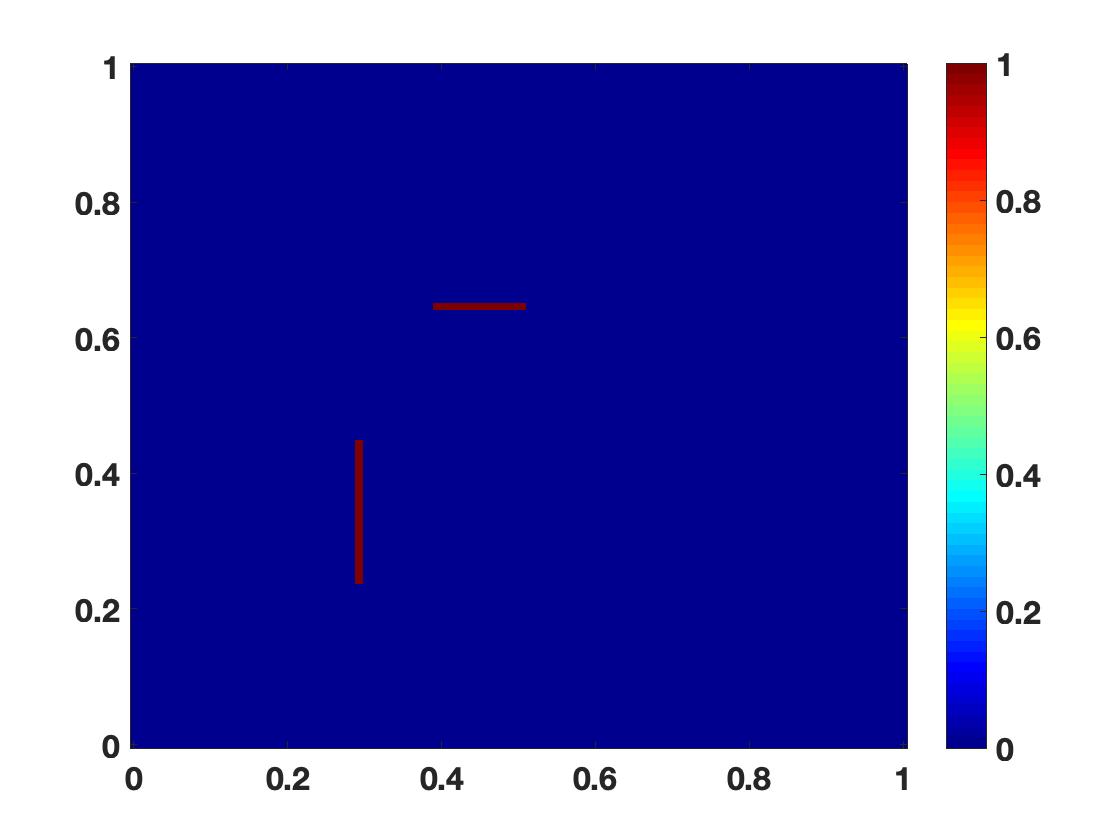}
		\caption{The heterogeneous permeability field: $\kappa$ (left) and source term $f$ (right).}
		\label{fig:kappa_2}
	\end{figure}

\begin{figure}[H]
		\centering
		\includegraphics[trim={2.5cm 1.5cm 1.5cm 0.0cm},clip,width=0.28 \textwidth]{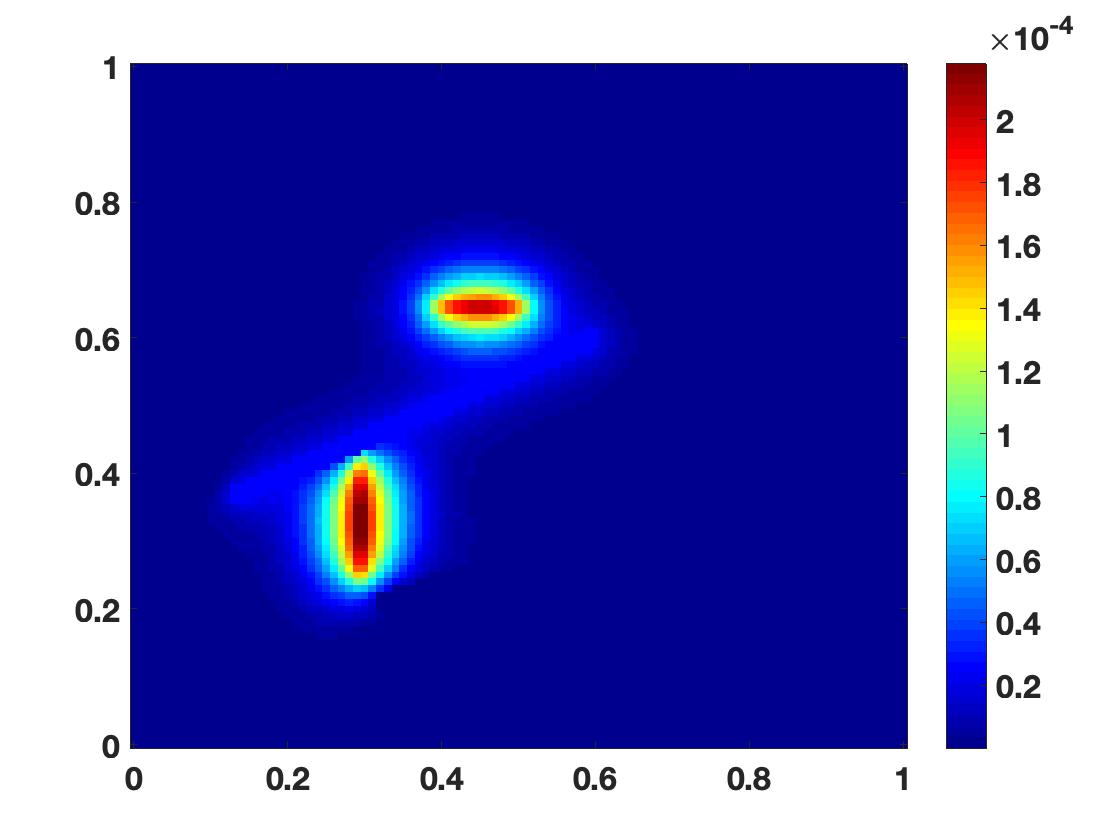}
		\includegraphics[trim={2.5cm 1.5cm 1.5cm 0.0cm},clip,width=0.28 \textwidth]{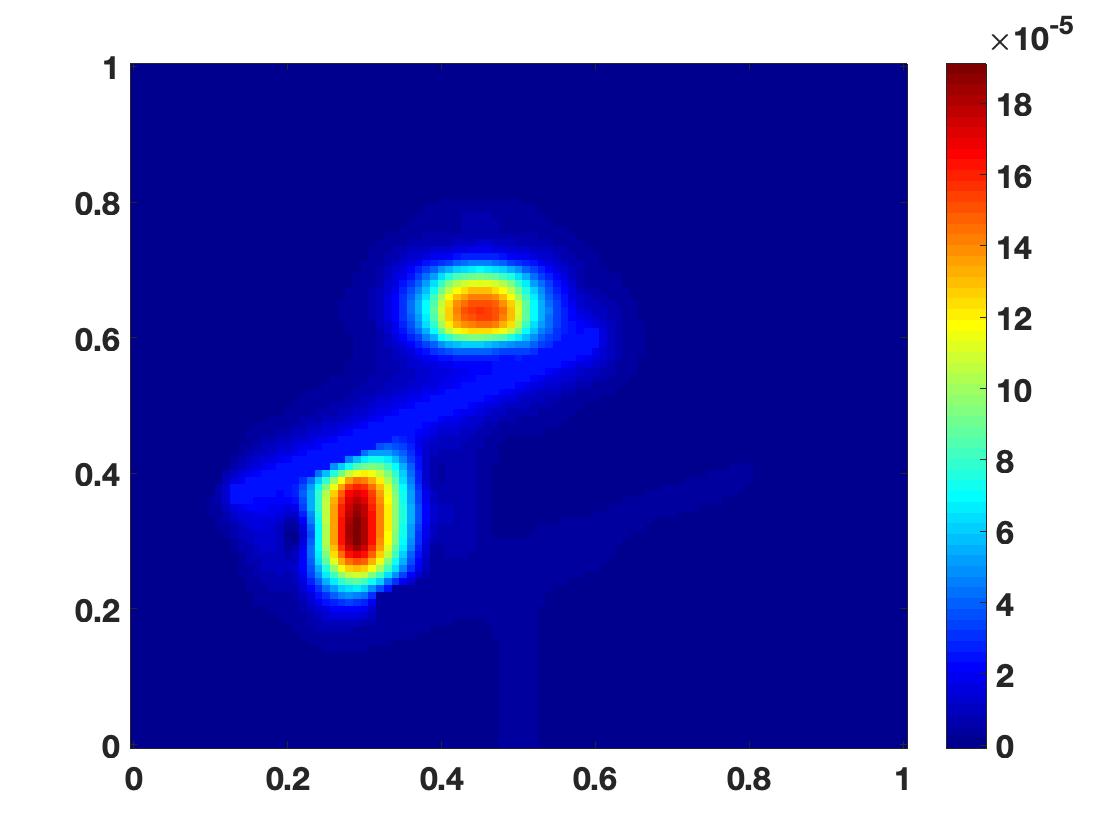}\\
		\includegraphics[trim={2.5cm 1.5cm 1.5cm 0.0cm},clip,width=0.28 \textwidth]{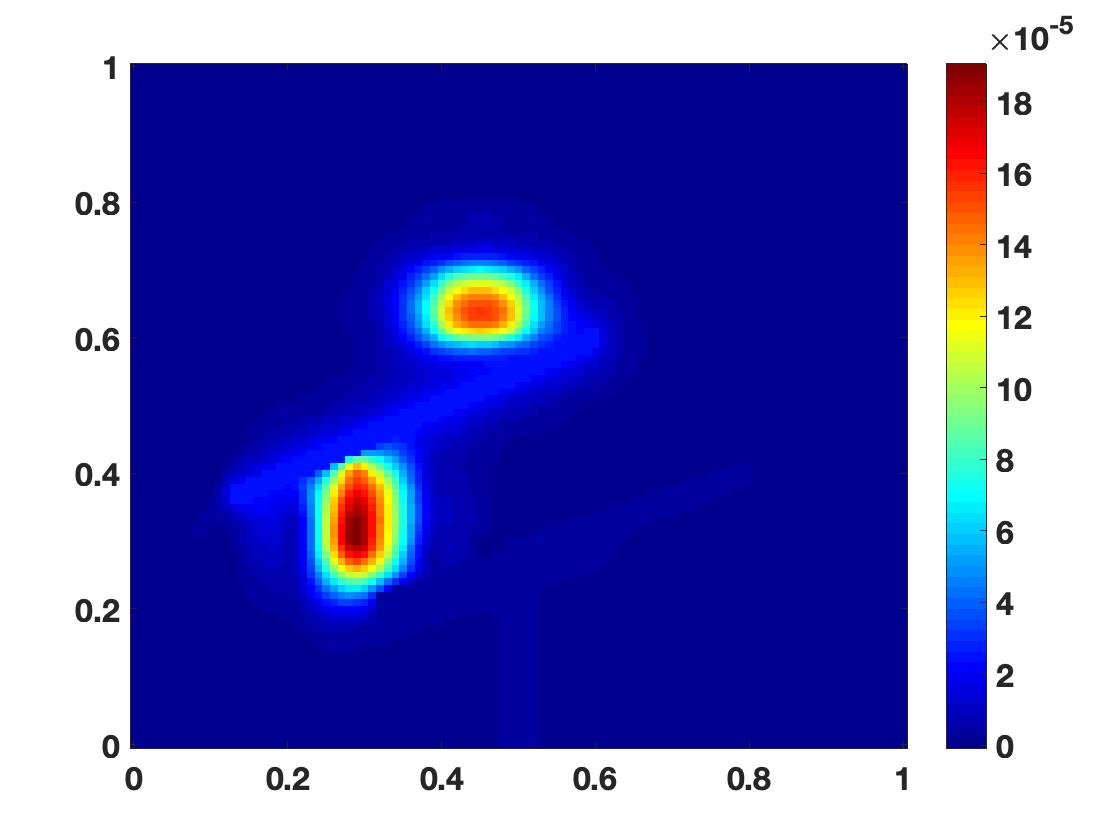}
		\includegraphics[trim={2.5cm 1.5cm 1.5cm 0.0cm},clip,width=0.28 \textwidth]{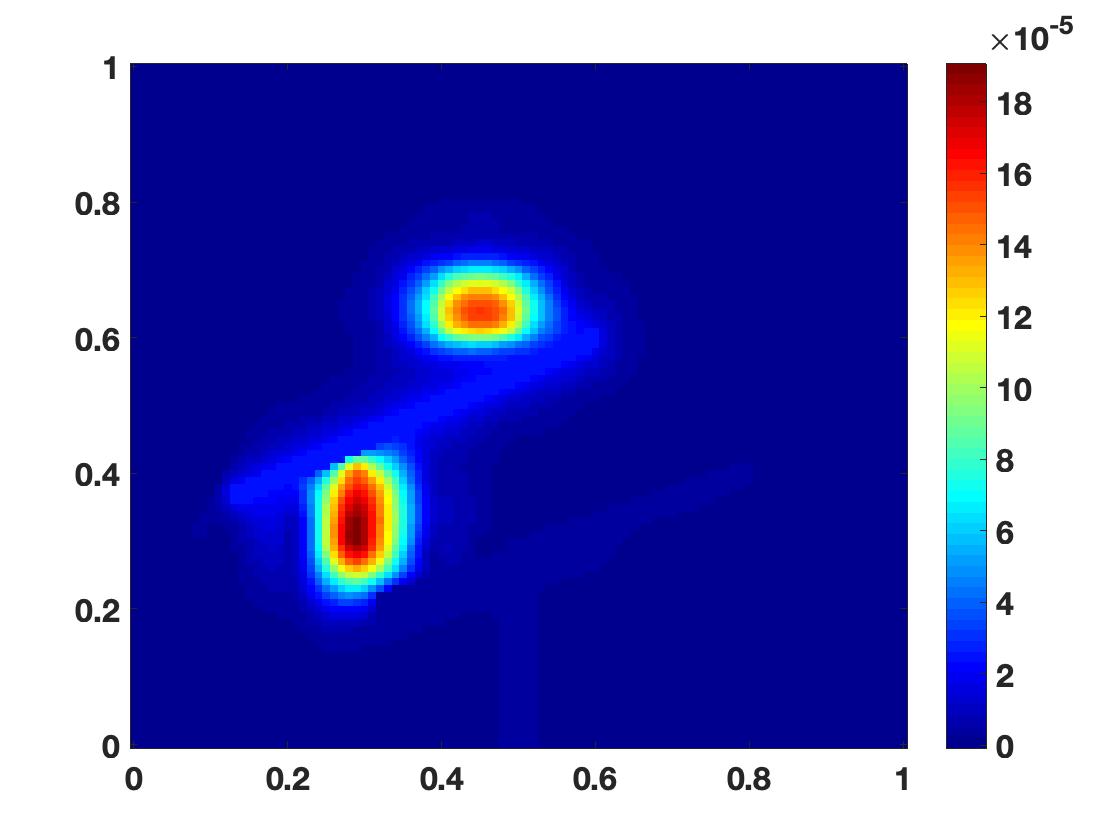}
		\caption{fine-grid solution $u_{h}^{2500}$ (top left), CEM solution $U_{cem}^{500}$ (top right), CEM solution with additional basis functions implicitly $\tilde{U}^{500}$ (bottom left), SCEM solution with additional basis functions explicitly $u_{scem}^{500}$  (bottom right).}
		\label{fig:solution_test2}
	\end{figure}

In the second experiment, we choose a heterogeneous permeability coefficient $\kappa$, which has two distinct value: 1 and $10^5$.
The source term $f$ is chosen to be a discontinuous function.
 They are plotted in Figure \ref{fig:kappa_2} for an illustration.
We first choose The fractional derivative order $\alpha:=0.9$. The fine-grid solution, CEM solution, CEM solution with more basis functions and SCEM solution at the final time $T$ are plotted in Figure \ref{fig:solution_test2}.

\begin{figure}[H]
		\centering
		\includegraphics[trim={2.5cm 0.0cm 1.0cm 0.2cm},clip,width=0.35 \textwidth]{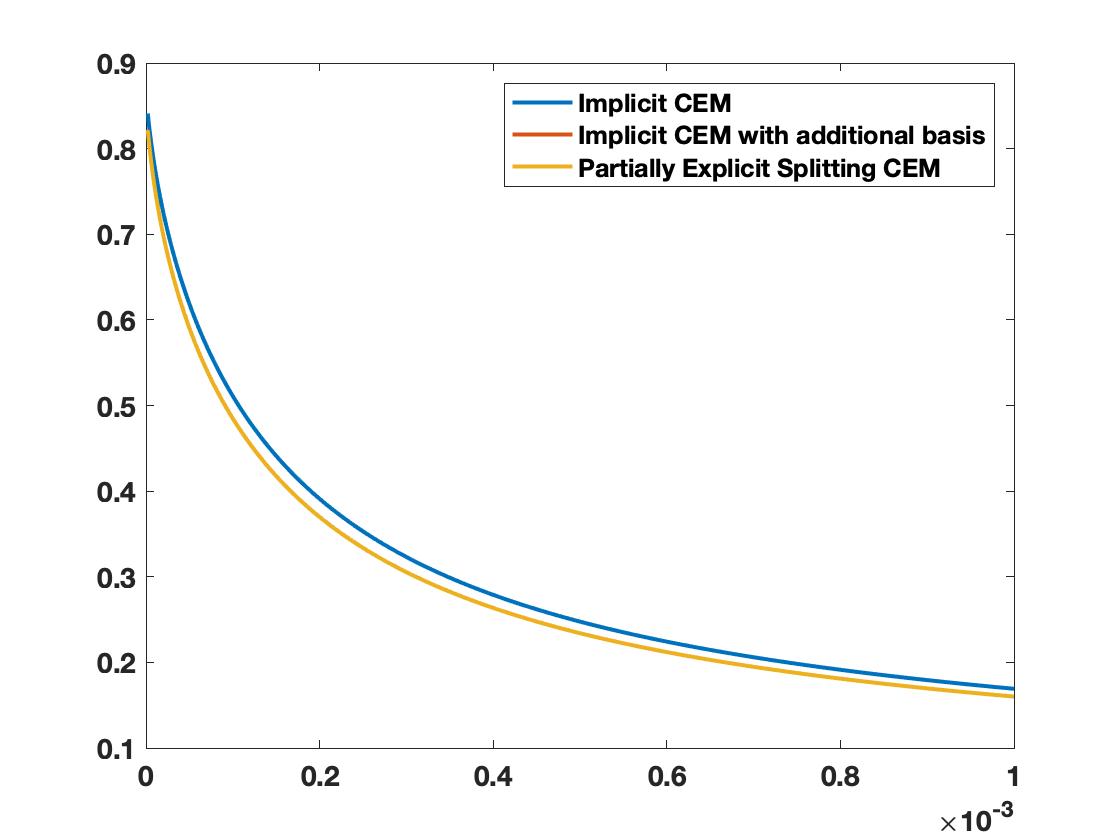}
		\includegraphics[trim={2.5cm 0.0cm 1.0cm 0.2cm},clip,width=0.35 \textwidth]{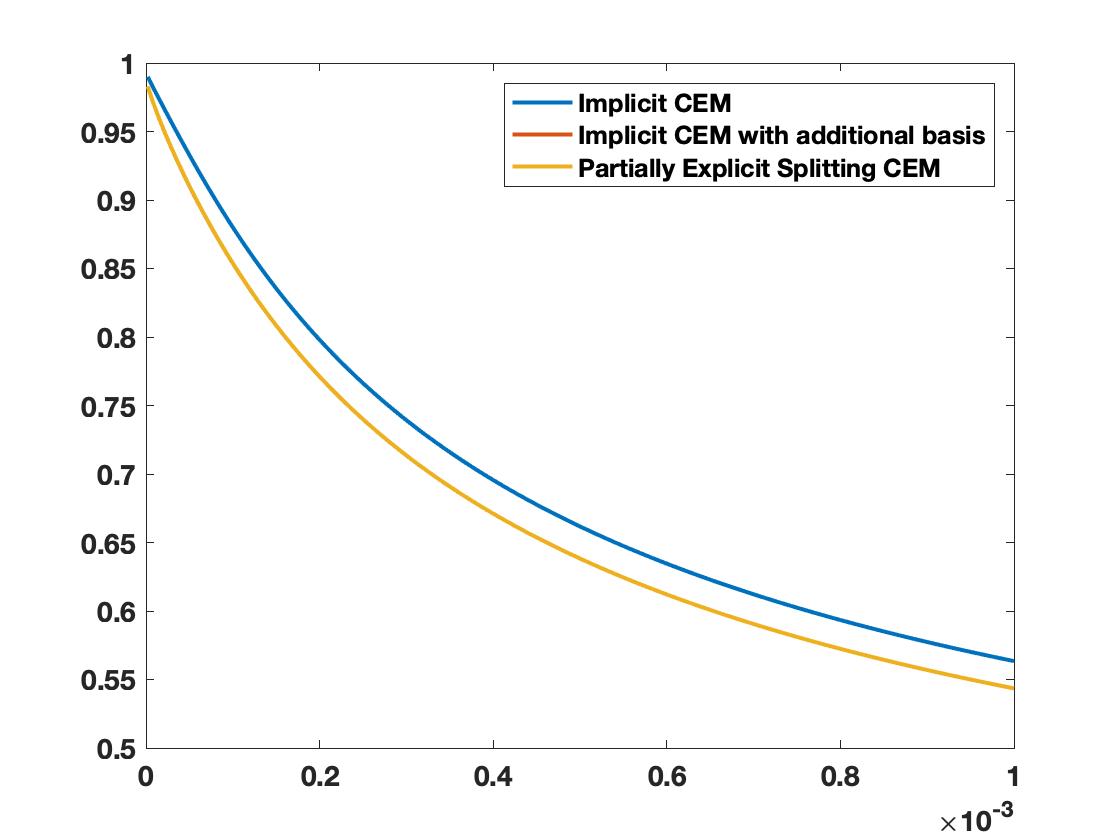}
		\caption{L2 error (left) and Energy error (right).}
		\label{fig:error_test2}
	\end{figure}

The convergence history of three numerical solutions in relative $L^2(\Omega)$-norm and relative $H^1_{\kappa}(\Omega)$-norm are presented 
 in Figure \ref{fig:error_test2}. In this experiment, we have similar observations as Experiment 1. That is, $u_{scem}$ has about similar accuracy as $\tilde{U}$ and it is cheaper to compute $u_{scem}$ than $\tilde{U}$.

\section{Conclusions}\label{sec:conclusions}

In this paper, we present a framework for partial explicit discretization
for time fractional PDEs. The work is motivated by many applications of 
multiscale time fractional PDEs. 
Explicit time discretization for time fractional PDEs requires very small time
steps due to size of the fine grid, the contrast,
 and additional power that is associated
with the fractional power of the time derivative.
Our approach solves time fractional PDEs
on a coarse grid by constructing appropriate coarse spaces.
We show that the proposed method is stable and one can choose the time
step that does not depend on the contrast and only depends on the coarse
mesh size. We note that our approach does not remove the constraint 
related to the power of the time fractional PDE.
 Via the construction of appropriate spaces
and careful stability analysis, we can show that the time step can 
be chosen not to depend on the contrast and scale as the coarse mesh size.
We present numerical results by considering time fractional diffusion 
in highly heterogeneous media. We show that the proposed partial explicit
approach provides similar results compared to the fully implicit method,
where all degrees of freedom are treated implicitly.

	\bibliographystyle{abbrv}
	\bibliography{ref,references,references4,references1,references2,references3,decSol}

\end{document}